\documentclass[11pt]{amsart}
\usepackage{latexsym}
\usepackage{amsfonts}
\usepackage{amsthm}
\usepackage{graphicx}
\usepackage{amssymb}
\usepackage{amsmath}
\usepackage{amssymb}
\makeatletter

\author[]{Alexei G. Myasnikov}
\address{Department of Mathematics and Statistics,  McGill
University, Montreal QC H3A 2K6 Canada} \email{amiasnikov@gmail.com}

\author[]{Vladimir Shpilrain}
\address{Department of Mathematics, The City  College  of New York, New York,
NY 10031} \email{shpil@groups.sci.ccny.cuny.edu}

\thanks{Research of both authors was partially supported by
the NSF grant DMS-0405105.}

\begin{document}

\title[]{Some metric properties of automorphisms of  groups}

\medskip

\begin{abstract}
\noindent    Study of the dynamics of automorphisms
 of a   group is usually focused on their {\it growth}
and/or  finite orbits, including fixed points.
In this paper, we introduce properties of a different kind;
using somewhat informal language, we call them {\it metric properties}.
 Two  principal  characteristics
of this kind are called here the ``curl"  and  the ``flux"; there seems
to be very little correlation between these and the growth
of an automorphism, which means they are likely to be an
essentially new tool for studying  automorphisms.

 We also observe that our definitions of the curl and  flux are
sufficiently general to be applied to mappings of arbitrary
metric spaces.
\end{abstract}

\maketitle

\section{Introduction}

 Let $G$ be a finitely generated  group of   rank
$r \ge 2$  with a  set  $X = \{x_1, ..., x_r \}$  of   generators,
and let $|w|$ be the usual lexicographic length of an element $w \in G$
with respect to $X$.

Let $\varphi$ be an automorphism  (or, more generally, an endomorphism)
of $G$ that takes $x_i$ to  $y_i$, $i= 1, ..., r$.
The  {\it growth function} of $\varphi$ with respect to $X$
can be   defined as

$$\Gamma_{\varphi, m}(n) = \max_{|w|=m} |\varphi^n(w)|.$$

 This function therefore measures, to some extent,  how fast the
length of elements of  $G$ can possibly increase under repeated
action of $\varphi$.

 One can also define a cumulative characteristic, usually called the
{\it growth rate}, or simply {\it growth}, of $\varphi$:

$$\Gamma(\varphi) = \sup_m \limsup_{n \to \infty} \root{n}\of{\Gamma_{\varphi, m}(n)}.$$

 For known properties of  growth of automorphisms
of a {\it free group} we refer to \cite{BFH}, \cite{Dicks}, and  \cite{LL}.
Very little seems to be known if  $G$ is not a free group.
\smallskip

In this paper, we introduce essentially new  characteristics
of an automorphism. These  will tell us how ``active" an automorphism is
rather
than how  it ``grows".

\medskip

\noindent  {\bf (1)} {\it Curl function} is  defined as

$$Curl_\varphi(n) = | \varphi(B_n) \cap B_n | =
\#\{w \in G, ~|w| \le n, ~|\varphi(w)| \le n\},$$

\noindent where $B_n$ is the ball of radius $n$ in the Cayley graph of
$G$.  This function therefore counts the number of elements left inside the ball
of radius $n$ by the automorphism $\varphi$.

As with the growth rate, one can define the  ``curl rate", or simply
``curl", of $\varphi$ as

$$Curl(\varphi) = \limsup_{n \to \infty} \root{n}\of{\frac{Curl_\varphi(n)}{|B_n|}}.$$

\medskip

\noindent  {\bf (2)} {\it Flux function} of $\varphi$ is  defined as

$$Flux_\varphi(n) = |B_n \setminus (\varphi(B_n) \cap B_n)| =
\#\{w \in G, ~|w| \le n, ~|\varphi(w)| > n\}.$$

 This function therefore counts the number of elements taken out of
  the ball $B_n$ of radius $n$ by the automorphism $\varphi$.

   Again, one can define the  ``flux rate", or simply ``flux", of $\varphi$
as follows:

$$Flux(\varphi) = \limsup_{n \to \infty} \root{n}\of{\frac{Flux_\varphi(n)}{|B_n|}}.$$

\smallskip

  We note that all these concepts can be defined for arbitrary {\it endomorphisms},
not necessarily automorphisms.
\smallskip

 It is immediately obvious that:
\medskip

\noindent  {\bf (i)} $0 \le Curl(\varphi), ~Flux(\varphi) \le 1$.
It is a very interesting question what values $Curl(\varphi)$ and $Flux(\varphi)$
can actually take. In Section \ref{props}, we show that there are gaps on the
scale of these values; in particular, $Flux(\varphi)$ cannot take values strictly
between 0 and $\frac{1}{4}$ for any injective endomorphism $\varphi$ of $F_r$.
\medskip

\noindent  {\bf (ii)} For any $n$,  $Curl_\varphi(n) + Flux_\varphi(n) = |B_n|$,
 the cardinality of the  ball $B_n$.
However, $Curl(\varphi) + Flux(\varphi) \ne 1$ in general; we shall see
relevant examples (e.g. Example 3.2) in Section \ref{examples}.
\medskip

 There are other, less obvious, properties of curl and flux  that we have collected
in Section \ref{props}. Whenever we give a particular property, we use it to
compare curl and flux  to growth. As it turns out, curl and flux have some
useful properties that growth does not have. For instance, we have
$Flux(\varphi) = Flux(\varphi^{-1})$ and $Curl(\varphi) = Curl(\varphi^{-1})$
for any automorphism $\varphi$ (Proposition 4.6); we also have some
inequalities for curl and flux functions of composite endomorphisms, including
$Flux_{\alpha \beta}(n) \le Flux_{\alpha}(n) + Flux_{\beta}(n)$
(Proposition 4.6), etc.

 We note at this point that Kaimanovich, Kapovich, and  Schupp \cite{KS}
have  independently come up with yet another dynamical characteristic of
an automorphism; they call it the {\it generic stretching factor}.
This is a number $\lambda = \lambda(\varphi)$ such that a given
automorphism $\varphi$ ``stretches" the length of ``almost all"
elements of  the  group approximately by a factor of   $\lambda$
(for more details see our Section \ref{estim}).
This stretching factor appears to be related (although not  directly)
 to our flux. In particular,  it is shown in \cite{KS} that
 the flux of any automorphism $\varphi$ of a {\it free group} is 1,
unless  $\varphi$ is a permutation of the set $X \cup X^{-1}$.
Moreover, if $\varphi$ is not a composition of
an inner automorphism and a permutation of the set $X \cup X^{-1}$,
then $\lim_{n \to \infty} \frac{Flux_\varphi(n)}{|B_n|}=1$.
Therefore, the flux cannot be used to distinguish
automorphisms of a free group.

 The situation with the curl however is different.   We show, for example,
that if  $Curl(\varphi)=1$, then  $\varphi$ is a composition of
an inner automorphism and a permutation of the set $X \cup X^{-1}$
(Theorem 5.1 in Section \ref{estim}). We also show that
``stabilizing"
 an automorphism of a free group (by  expanding the free
generating set $X$) may  change its  curl, but not the  growth
(see  Example 3.4 in Section \ref{examples}).  This is, arguably, an evidence
of the  curl being a
more delicate  characteristic of an automorphism than its growth.

  To conclude the Introduction, we observe that our definitions of curl
and flux are sufficiently general to be applied to mappings of
arbitrary metric spaces.

\section{Problems}
\label{problems}

In this section, we  list a few open problems that are,
 in our opinion, important for better understanding the nature of
curl and flux. As usual,  $F_r$ denotes the free group of rank
$r \ge 2$ with a set  $X$ of free generators.
\medskip

\noindent  {\bf Problem 1.} {\bf (a)} What is the maximum (or supremum) of
possible $\ne 1$ values
of the curl for  automorphisms (endomorphisms)   of $F_r$ ?
\smallskip

\noindent  {\bf (b)} What is the minimum (or infinimum)
 of possible values  of the curl for  automorphisms   of $F_r$ ?

\medskip

 A good start would be $r=2$. It is conceivable that the automorphism
$\alpha: x \to xy, ~y \to y$ has the maximum possible $\ne 1$ curl
among automorphisms of $F_2$, but we do not have a proof of that.
Nor  do we have the exact value of $Curl(\alpha)$; according to computer
experiments (see Section \ref{estim}), this value is approximately 0.956.

  We also note here that the infinimum of possible values  of the curl for
{\it endomorphisms} of $F_r$ is $\frac{1}{2r-1}$, see Proposition 4.1
 in Section \ref{props}.
\medskip

\noindent  {\bf Problem 2.} What is the minimum (or infinimum)
 of possible positive values
of the flux for   endomorphisms   of $F_r$ ?
\medskip

 As we have mentioned in the Introduction, $Flux(\varphi)$ cannot take values
strictly between 0 and $\frac{1}{4}$ for any injective  endomorphism of $F_r$.
If $\varphi$ is an automorphism of $F_r$, then $Flux(\varphi)=0$ or $1$ by the
result of \cite{KS} mentioned in the Introduction.  This is however not
the case for arbitrary
endomorphisms; for example, the endomorphism of $F_2$ given by
$x \to xy, ~y \to 1$ has the flux strictly between 0 and  1 (see Example 3.7
 in Section \ref{examples}).

\medskip

\noindent  {\bf Problem 3.} Are values of flux and curl always
algebraic numbers? (Values of growth are.)
\medskip

\noindent  {\bf Problem 4.} Find the {\it exact} value of $Curl(\varphi)$
 for at least one  $\varphi \in Aut(F_r)$ with $Curl(\varphi) \ne  1$.
\medskip

\noindent  {\bf Problem 5.} Suppose $Curl(\varphi)=Curl(\psi)$ for some
automorphisms $\varphi, \psi$ of $F_r$. Is it true that $\varphi$ is
a composition of $\psi$ with a permutation of the set $X \cup X^{-1}$
and an inner automorphism?
\medskip

 The converse is true (see  Proposition 4.2 in Section \ref{props}).
If the answer to Problem 5 is affirmative,
this  will mean that the curl is indeed a very sharp  characteristic of
a free group automorphism. We were able to show that if
 $Curl(\varphi)=1$, then  $\varphi$ is a composition of
an inner automorphism and a permutation of the set $X \cup X^{-1}$
(Theorem 5.1 in Section \ref{estim}).

 The following problem is rather vague, but it appears to be important.
\medskip

\noindent  {\bf Problem 6.}  Find tight bounds for $Curl(\alpha \beta)$
in terms of $Curl(\alpha)$,  $Curl(\beta)$. More generally, what
information about $Curl(\alpha \beta)$ can be extracted from knowing
$Curl(\alpha)$ and  $Curl(\beta)$ ?

\section{Examples}
\label{examples}

 In this section, we compute curl and flux for some simple automorphisms
of $F_r$,  the free group of rank  $r \ge 2$ with a set  $X$ of
free generators.

\medskip

\noindent  {\bf Example 3.1.}  Let $\pi$ be any  automorphism
that   permutes the elements of the set $X \cup X^{-1}$.
 Then, since $\pi$ does not change the
length of any element, we have  $Curl(\pi)=1, ~Flux(\pi)=0$.
It is also obvious that the growth function of $\pi$ is identically
equal to 1. $\Box$
\medskip

\noindent  {\bf Example 3.2.} Let $i_g$ be the conjugation by
an element $g \in F_r$. Then $Curl(i_g)=Flux(i_g)=1$.
Indeed, it is sufficient to limit considerations to elements
of a sphere $S_n$  because these comprise ``most" of the
elements of the ball $B_n$ (see \cite{KMSS} for more rigorous
estimates supporting this claim).
Now suppose $g$ ends with $x$ for some
$x \in X \cup X^{-1}$. Then an element $u \in S_n$ gets
taken out of $B_n$ by $i_g$  if $u$ does not start with $x^{-1}$.
The number of elements with this property has the same growth
function, up to a constant factor, as the total number of
elements in $S_n$ does. This yields $Flux(i_g)=1$.

 On the other hand, an element $u \in S_n$ is {\it not}  taken out
of $B_n$ by $i_g$  if $u$ starts with $g^{-1}$. Again, the
number of elements with this property has the same growth
function, up to a constant factor, as the total number of
elements in $S_n$ does. This yields $Curl(i_g)=1$. $\Box$
\medskip

\noindent  {\bf Example 3.3.} Let $r=2$, and denote the generators of the
group $F_2$ by $x$ and $y$. Let $\alpha: x \to xy, ~y \to y$.
Then the growth function of $\alpha$ is easily seen to be linear in $n$,
whereas both $Curl_{\alpha}(n)$ and $Flux_{\alpha}(n)$ are exponential. $\Box$

\medskip

\noindent  {\bf Example 3.4.} Again, let $r=2$, and let $i_x$ be the
conjugation by the generator $x$. Then $Curl(i_x)=Flux(i_x)=1$. Now extend $i_x$
to the free group $F_3$ generated by $x$, $y$, and $z$, by fixing the
extra generator $z$. Call this new automorphism $\widehat{i_x}$. Thus,
$\widehat{i_x}: x \to x, ~y \to xyx^{-1}, ~z \to z$.
 Then, since $\widehat{i_x}$ is not a composition of
an inner automorphism and a permutation of the set $X \cup X^{-1}$, we have
$Curl(\widehat{i_x}) < 1$ by Theorem 4.1 in our Section \ref{estim}.
 $\Box$

\medskip

 Thus, Example 3.4 shows that the curl of an automorphism can change (decrease)
 under ``stabilization". This makes contrast with the growth and reinforces
 the impression that the curl  reflects more delicate
  properties of automorphisms than the growth does.

 In the next example, we show that the curl of an endomorphism can also
{\it increase} under ``stabilization".
\medskip

\noindent  {\bf Example 3.5.}  Let $r=2$, and let
$\varphi: x \to x^{5}, y \to y^{5}$
be  an endomorphism of the group $F_2$. Then, by  Proposition 4.1 in
Section \ref{props}, $Curl(\varphi) = \frac{3^{\frac{1}{5}}}{3}$.

 For computational convenience, let us now ``stabilize" $\varphi$ by
adding two extra generators,  $z$ and  $t$. Thus,
$\widehat{\varphi}: x \to x^{5}, y \to y^{5}, ~z \to z, ~t \to t$.
Then, for any $u=u(z, t)$ of length $n$, we have  $|\varphi(u)| = n$.
There are   at least  $3^n$ words $u$ like that.
Therefore, $Curl(\varphi) \ge \frac{3}{7} > \frac{3^{\frac{1}{5}}}{3}$.  $\Box$

\medskip

\noindent  {\bf Example 3.6.}  Again, let $r=2$, and let
 $\varphi = \alpha \cdot \pi_{xy}$, where  $\alpha: x \to xy, ~y \to y$, and
$\pi_{xy}$ permutes $x$ and $y$. Thus,  $\varphi: x \to xy, ~y \to x$.
 Then it is fairly clear that $\varphi$ has exponential growth
(i.e.,  $\Gamma(\varphi) > 1$), whereas
$\alpha$ has linear growth (in particular,  $\Gamma(\alpha) = 1$).
At the same time, $Flux(\varphi) = Flux(\alpha)$
and $Curl(\varphi) = Curl(\alpha)$ since $\varphi$ is a composition
of $\alpha$ with a length-preserving automorphism. $\Box$

  The point of this example is to show, again, that the curl and
the flux of an automorphism seem to have very little or no correlation
with the growth.
\medskip

 We conclude this section with an example of an endomorphism
$\varphi$ of the group $F_2$ whose flux is strictly between 0
and  1.

\medskip

\noindent  {\bf Example 3.7.} Let $\varphi: x \to xy, ~y \to 1$.
Then  $0 < Flux(\varphi) < 1$.  Indeed, if a word $w$ of length $n$
has $> \frac{n}{2}$ occurrences of $x$ and no occurrences of $x^{-1}$,
then $|\varphi(w)|>n$. The number of words like that is at least
 ${n}\choose{\frac{n}{2}}$, which is exponential in  $n$.
This shows that $0 < Flux(\varphi)$.

 To show $Flux(\varphi) < 1$, we observe that for a word $w$ of length $n$
to be taken out of $B_n$ by $\varphi$, it should have the exponent sum
on $x$ greater than  $\frac{n}{2}$ (by the absolute value).
This implies that  the number of occurrences in $w$ of either $x^{-1}$
or $x$ should be $\ge \frac{3n}{4}$. The set of words like that is
exponentially negligible in $B_n$ by \cite[Proposition 6.1]{KSS}.
Therefore, $Flux(\varphi) < 1$. $\Box$

\section{Some properties of curl and flux}
\label{props}

 In this section, we gather some interesting, in our opinion,
properties of curl and flux.  Most of these properties are valid
for arbitrary {\it endomorphisms}, not necessarily automorphisms.

\medskip

\noindent  {\bf Proposition 4.1.} {\bf (a)}  Let $k \ge 2$, and  let
 $\varphi: x_i \to x_i^k, ~i=1, ..., r$
be  an endomorphism of the group $F_r$. Then
$Curl(\varphi) = \frac{(2r-1)^{\frac{1}{k}}}{2r-1}$.
\smallskip

\noindent  {\bf (b)} For any endomorphism $\psi$ of the group $F_r$,
$Curl(\psi) \ge \frac{(2r-1)^{\frac{1}{k}}}{2r-1}$ for some $k \ge 2$. Therefore,
the infinimum of possible values  of the curl for
 endomorphisms of $F_r$ is $\frac{1}{2r-1}$.

\medskip

\noindent  {\bf Proof.} {\bf (a)} Note that for any   $u \in F_r$,
one has $|\varphi(u)| = k |u|$. Therefore, $Curl_\varphi(n)$ is just
equal to the number of elements of length $ \le \frac{n}{k}$ in
$F_r$, i.e., to $O((2r-1)^{\frac{n}{k}})$, whence the result.
\smallskip

\noindent  {\bf (b)}  Let $\psi: x_i \to y_i, ~i=1, ..., r$, and
suppose  $|y_i| \le k$ for some  $k \ge 2$. Then $|\psi(u)|  \le
k|u|$ for any   $u \in F_r$. Therefore, whenever $|u| \le
\frac{n}{k}$, one has $|\psi(u)| \le n$. The result follows. $\Box$

\medskip

\noindent  {\bf Proposition 4.2.}
{\bf (a)}  Composing any endomorphism $\varphi$ of $F_r$
with any permutation of the set $X \cup X^{-1}$
does not change either $Flux(\varphi)$  or  $Curl(\varphi)$.
\smallskip

\noindent  {\bf (b)} Composing any endomorphism $\varphi$ of $F_r$
with any inner automorphism does not change  $Curl(\varphi)$.
If $\varphi$ is injective, then such composing
does not change $Flux(\varphi)$ either.
\medskip

\noindent  {\bf Proof.}  Part (a) is obvious, so we proceed with
part (b). Note that $Curl(\varphi) >0$ by Proposition 4.1 and
$Flux(\varphi) >0$ by Theorem 4.4 below. Then the argument similar
to that in Example 3.2 shows that if we apply $\varphi$ followed by
an inner automorphism, this will not change either
$Flux(\varphi)$  or  $Curl(\varphi)$.

 Suppose now an inner automorphism is applied first, followed by
$\varphi$. By using inductive argument, we may assume, to simplify
the notation, that the inner automorphism is $i_x$, i.e. conjugation
by $x \in X$. Then $i_x$ leaves inside $B_n$  all elements $v \in B_n$
that start with $x^{-1}$. Suppose now  an element $w \in B_n$
starts with some other $y \in X\cup X^{-1}$, i.e., $w = yu$. If this
$w$ is left inside $B_n$ by $\varphi$, then so is $w^{-1} = u^{-1}y^{-1}$.
The number of elements in $B_n$ of the form $u^{-1}y^{-1}$ is
the same, up to a constant factor, as the number of elements of the form
$x^{-1}u$. Each of these numbers is equal, again up to a constant factor,
to the total number of elements in $B_n$. These two facts show that
the curl of the composite endomorphism is the same as the curl of $\varphi$.
 The flux is treated similarly. $\Box$

\medskip

 Before we get to the next  result, we need a lemma:
\medskip

\noindent  {\bf Lemma 4.3.} Let $\varphi$ be an endomorphism of $F_r$
such that,
for some cyclically reduced $v \in F_r$, one has $|\varphi(v)| \ge 2|v|$.
Then $Flux(\varphi) \ge \frac{1}{4}$.
\medskip

\noindent  {\bf Proof.}  By Example 3.2, we may assume that $\varphi$ is
{\it not} a conjugation.
We are going to fix a particular  $k$ and build
sufficiently many words  $w \in F_r$ of length $k$ whose length is
increased by $\varphi$. To that effect, we first fill in the
leftmost  $\ge \frac{k}{2}$ positions with $v^s$, where
$s= [\log_{|v|}\frac{k}{2}]+1$. Let $m=|v^s|-\frac{k}{2}$; then
$0 \le m \le  |v|$.

 Now we designate the rightmost
$\frac{k}{4} - \frac{m}{2} - 1$
positions in $w$ as ``arbitrary" (call this part $w_{right}$),
 and fill in the intermediate $\frac{k}{4} -\frac{m}{2} + 1$ positions
as follows:
\medskip

\noindent  {\bf (i)} Among all words in $F_r$ of length
$\frac{k}{4} -\frac{m}{2} - 1$
choose one, call it $u$,  such that $|\varphi(u)| \ge |\varphi(g)|$ for any
$g$ of length $\frac{k}{4} -\frac{m}{2} - 1$,  and place $u$ immediately
left of $w_{right}$.
 That way, after we apply $\varphi$ to $w$,
cancellation between $\varphi(w_{right})$  and  $\varphi(u)$ cannot possibly go
left beyond $\varphi(u)$.
\medskip

\noindent  {\bf (ii)}  Fill in the remaining two positions right of the
 $v^s$ with two letters, call them  $a$  and  $b$,
in such a way that there is no cancellation between either
$\varphi(v^s)$ and $\varphi(ab)$, or  between $\varphi(ab)$
and $\varphi(u)$,  or  between $\varphi(a)$ and $\varphi(b)$
(this is possible since $\varphi$ is  not a conjugation).
Then the length of $\varphi(w)$ is greater than $k$.

 Finally, we observe that the number of different $w_{right}$ of length
$\frac{k}{4} -\frac{m}{2} - 1$ grows as $r^\frac{k}{4}$, up to an
exponentially negligible factor.
This yields the result. $\Box$

\medskip

\noindent  {\bf Theorem 4.4.} Let $\varphi$ be
an injective endomorphism  of $F_r$. Then
either

\noindent  {\bf (a)} $Flux(\varphi) =0$, in which case  $|\varphi(x)| = 1$
 for all $x \in X$,

or

\noindent  {\bf (b)} $Flux(\varphi) \ge \frac{1}{4}$.
\medskip

\noindent  {\bf Proof.} If $|\varphi(x)| = 1$
 for all $x \in X$, then obviously $Flux(\varphi) =0$.
Let now $\varphi(x_i)=y_i, ~|y_i| \ge 2$ for some $i$. Consider two cases:
\smallskip

\noindent  {\bf (1)} For some $i$, $|y_i| \ge 2$ and  $y_i$ is cyclically reduced.
Then  $Flux(\varphi) \ge \frac{1}{4}$ by Lemma 4.3 if we let $v = x_i$.

\smallskip

\noindent  {\bf (2)} For all $i$ such that $|y_i| \ge 2$,
one has  $y_i$ not cyclically reduced. Here we have two subcases:

\smallskip

\noindent  {\bf (i)} there are $k, l, ~k \ne l$, such that for some
$i, j$ one has $y_i=x_kg_ix_k^{-1}, ~y_j=x_l g_j x_l^{-1},$ and
at least one of the $y_i, y_j$ has length $\ge 2$.  Then, for $u = x_i x_j$,
we have $|\varphi(u)| \ge 2|u|$
 and $u$ is cyclically reduced. Then, by Lemma 4.3, we have
$Flux(\varphi) \ge \frac{1}{4}$.
\smallskip

\noindent  {\bf (ii)} every $y_i$ with $|y_i| \ge 2$ is of the form
$xg_ix^{-1}$ for some fixed $x \in X \cup X^{-1}$. Suppose, for some $j$,
$\varphi(x_j)=x_k \ne x$. Then, for $u = x_i x_k$, we have $|\varphi(u)| \ge 2|u|$
 and $u$ is cyclically reduced. Then, by Lemma 4.3, we have
$Flux(\varphi) \ge \frac{1}{4}$. The remaining case is where {\it every}
$y_i$ is of the form $xg_ix^{-1}$. If, for some $i$, $|g_i| \ge 2$, then
the argument from the proof of Lemma 4.3 will
still work after obvious minor adjustments. If $|g_i| = 1$ for every $i$,
then $\varphi$ is a composition of a permutation with the conjugation by
$x$, whence $Flux(\varphi)=1$. $\Box$

\medskip

\noindent  {\bf Proposition 4.5.} For any automorphism $\varphi$
of any group $G$,
$Flux(\varphi) = Flux(\varphi^{-1})$ and   $Curl(\varphi) = Curl(\varphi^{-1})$.
Moreover, for any $n \ge 1$, $Flux_{\varphi^{-1}}(n) = Flux_{\varphi}(n)$,
and $Curl_{\varphi^{-1}}(n) = Curl_{\varphi}(n)$.
\medskip

\noindent  {\bf Proof.} Let $A$ be the set of elements of $B_n$ taken out of
$B_n$  by  $\varphi$, and $B$ the set of elements of $B_n$
left by $\varphi$ inside the ball. Furthermore, let $C$ be the set of elements
outside of $B_n$ taken by  $\varphi$ inside $B_n$,  $A'$  the set of elements
of $B_n$ taken out of $B_n$  by  $\varphi^{-1}$, and $B'$ the set of elements of $B_n$
left by $\varphi^{-1}$ inside the ball.

 Then, since $\varphi$ is {\it onto}, we must have
 $|C| = |A|$. At the same time, we clearly have
$|C| = |A'|$, hence $|A'| = |A|$. This implies $Flux_{\varphi^{-1}}(n) = Flux_{\varphi}(n)$.

 Now since $A \cup B = A' \cup B' = B_n$ and $A \cap B = A' \cap B' = \emptyset$,
we have $|B'| = |B|$, whence  $Curl_{\varphi}(n) = Curl_{\varphi^{-1}}(n)$. $\Box$

\medskip

\noindent  {\bf Proposition 4.6.} For any automorphisms $\alpha$ and  $\beta$
of any group $G$ and for any $n \ge 1$, one has:
\smallskip

\noindent  {\bf (a)}  $Curl_{\alpha \beta}(n) \le Curl_{\beta}(n) + Flux_{\alpha}(n)$.

\smallskip

\noindent  {\bf (b)}  $Flux_{\alpha \beta}(n) \le Flux_{\alpha}(n) + Flux_{\beta}(n)$.
\smallskip

\noindent  {\bf (c)}  $Curl_{\alpha \beta}(n) \ge Curl_{\beta}(n) - Flux_{\alpha}(n)$.
\smallskip

\noindent  {\bf (d)}  $Flux_{\alpha \beta}(n) \ge  Flux_{\beta}(n) - Flux_{\alpha}(n)$ and
 $Flux_{\alpha \beta}(n) \ge  Flux_{\alpha}(n) - Flux_{\beta}(n)$.
\smallskip

\noindent  {\bf (e)}  $Flux_{\alpha \beta}(n) \ge  Curl_{\alpha}(n) - Curl_{\beta}(n)$.
\smallskip

\noindent  Inequalities (a) and  (b) are actually valid for arbitrary
endomorphisms.

\medskip

\noindent  {\bf Proof.} First of all, we note that when we write $\alpha \beta$,
we assume that $\alpha$ is applied first.

\noindent  {\bf (a)}  Elements left inside the ball  $B_n$ by the
automorphism $\alpha \beta$ are among  those left inside $B_n$ by $\beta$
or  among  those first taken by $\alpha$ outside $B_n$, and then taken
back inside by $\beta$. The quantity of the  former is bounded by
$Curl_{\beta}(n)$, and the quantity of the latter by $Flux_{\alpha}(n)$.
This completes the proof of part (a).
\smallskip

\noindent  {\bf (b)} Argument similar to the one in (a) establishes
this inequality.
\smallskip

\noindent  {\bf (c)}  In (a), plug in $\alpha^{-1}$ for $\alpha$ and
$\alpha \beta$ for  $\beta$.  Then observe that
$Flux_{\alpha}(n) = Flux_{\alpha^{-1}}(n)$ by Proposition 4.5.
\smallskip

\noindent  {\bf (d)} Re-write (b) as
$Flux_{\alpha}(n) \ge  Flux_{\beta}(n) - Flux_{\alpha \beta}(n)$.
Now plug in $\alpha \beta$ for  $\alpha$ and  $\beta^{-1}$ for  $\beta$
to get $Flux_{\alpha \beta}(n) \ge Flux_{\beta^{-1}}(n) - Flux_{\alpha}(n)$.
Since, by Proposition 4.5, $Flux_{\beta^{-1}}(n) = Flux_{\beta}(n)$, this
yields the first inequality.

 For the second inequality, plug in $\alpha \beta$ for $\alpha$ and
$\beta^{-1}$ for  $\beta$ in (b). Then we get $Flux_{\alpha}(n) \le Flux_{\alpha \beta}(n)
+ Flux_{\beta}^{-1}(n)$. Since $Flux_{\beta^{-1}}(n) = Flux_{\beta}(n)$
by Proposition 4.5, this yields the result.
\smallskip

\noindent  {\bf (e)}  In (a), plug in $\alpha \beta$ for $\alpha$ and
$\beta^{-1}$ for  $\beta$. Then observe that
$Curl_{\beta}(n) = Curl_{\beta^{-1}}(n)$ by Proposition 4.5. $\Box$

%\noindent  {\bf Theorem 4.7.} Let $\varphi$ be an automorphism of $F_r$
%which is not a permutation of the set $X \cup X^{-1}$. Let $\widehat{\varphi}$ be
%the automorphism of $F_{r+1}$ obtained as follows: $\widehat{\varphi}(x_i)=
%\varphi(x_i), ~i=1, ..., r; ~\widehat{\varphi}(x_{r+1})=x_{r+1}.$
%Then $Curl(\widehat{\varphi}) < Curl(\varphi)$.

\section{Evaluating the curl}
\label{estim}

 Computing the exact value of   $Curl(\varphi)$ is  a
difficult problem for most automorphisms $\varphi$ of a free group,
so the best one can hope for (at least for now)  is to somehow estimate
 that  value. In this section, we are able to give the affirmative
answer to Problem 5 from Section \ref{problems} in the special case where
$\psi$ is the identity automorphism.

 To fully appreciate Theorem 5.1 below, the reader should bear in mind
that, according to computer
experiments (see the tables in the end of this section), for the
automorphism $\varphi: x \to xy, ~y \to y$ of $F_2$,
$Curl(\varphi)$ is approximately 0.956.

\medskip

\noindent  {\bf Theorem 5.1.}  Let $\alpha$ be an automorphism
of the group $F_r$ which is not a composition of
an inner automorphism and a permutation of the set $X \cup X^{-1}$.
 Then $Curl(\alpha) < 1$. Moreover,  $Curl(\alpha)$ is bounded away
 from 1, i.e., there is a   constant $c=c(r)<1$, independent of
 $\alpha$, such that $Curl(\alpha) < c$.
\medskip

\noindent  {\bf Proof. } To simplify the language, let us call
an automorphism {\it simple} if it is a composition of
an inner automorphism and a permutation of the set $X \cup X^{-1}$.

 Denote by $M_n(\alpha)$ the set $\{w \in F_r, ~|w| \le n, ~|\alpha(w)| \le n\}$.
 Recall that the cardinality of this set is what we call the curl function
$Curl_\alpha(n)$ of $\alpha$.

 Let $\lambda >1$. Clearly,
$$M_n(\alpha) ~\subseteq  ~B_{\frac{n}{\lambda}} \cup
\{u \in F_r, ~\frac{n}{\lambda} < |u| \le n, ~|\alpha(u)| \le \lambda |u|\}.$$

\noindent   The first set in the union on the right, the ball of
radius $\frac{n}{\lambda}$, is {\it asymptotically exponentially
negligible compared to $B_n$} (or just {\it asymptotically
exponentially negligible}, to simplify the language), which means

$$\lim_{n \to \infty} \root{n}\of{\frac{|B_{\frac{n}{\lambda}}|}{|B_n|}}<1.$$

 By   \cite[Theorem 6.8]{KS}, if  $\lambda <1 + \frac{2r-3}{2r^2-r}$,
then, since  $\alpha$ is not simple, the second  set in the union
above, i.e., the set
$$S_{\lambda, \alpha}(n)=\{u \in F_r, ~\frac{n}{\lambda}
< |u| \le n, ~|\alpha(u)| \le \lambda |u|\},$$

\noindent  must be asymptotically  exponentially negligible, too.

 Since  the  union of two asymptotically  exponentially negligible sets  is itself
asymptotically exponentially negligible, this implies that the set
$M_n$ is asymptotically  exponentially negligible, hence
$Curl(\alpha) < 1$.

To  prove the last claim in the statement of Theorem 5.1, we note
that, for a fixed  $\lambda$ such that $1< \lambda <1 +
\frac{2r-3}{2r^2-r}$, both the limits ~$\lim_{n \to \infty}
\root{n}\of{\frac{|B_{\frac{n}{\lambda}}|}{|B_n|}}$ ~and ~$\lim_{n
\to \infty} \root{n}\of{\frac{|S_{\lambda, \alpha}(n)|}{|B_n|}}$
~are bounded away  from 1 by a   constant $c=c(r)<1$, independent of
 $\alpha$. For the former limit, this is obvious. For the latter
 limit, this follows from the argument in the beginning of the proof
 of \cite[Theorem 6.8]{KS}.
$\Box$

\bigskip

 In conclusion, we present the results of computer experiments on evaluating
flux and curl of several automorphism. In the tables below, we give
values of the {\it curl ratio}
$\frac{Curl_\varphi(n)}{|B_n|}$ and
   the {\it flux ratio} $\frac{Flux_\varphi(n)}{|B_n|}$  along with
 the {\it curl root}
$\root{n}\of{\frac{Curl_\varphi(n)}{|B_n|}}$ and
   the {\it flux root} $\root{n}\of{\frac{Flux_\varphi(n)}{|B_n|}}$.

 We start with the ``simplest non-simple" automorphism of $F_2$.

$$
\left\{
\begin{array}{l}
x \rightarrow xy \\
y \rightarrow y
\end{array}
\right.
$$

\begin{center}
\begin{tabular}[]{r||l|l|l|l}
 \mbox{ n} & \mbox{\tiny CURL\_RATIO} & \mbox{\tiny CURL\_ROOT} & \mbox{\tiny FLUX\_RATIO} & \mbox{\tiny FLUX\_ROOT} \\
 \hline
   10 &     0.331634 &     0.895501 &     0.668366 &     0.960509 \\
   20 &     0.181176 &     0.918132 &     0.818824 &     0.990055 \\
   50 &    0.0372579 &      0.93632 &     0.962742 &     0.999241 \\
  100 &    0.0033803 &      0.94469 &      0.99662 &     0.999966 \\
  200 &  3.55979e-05 &     0.950073 &     0.999964 &            1 \\
  300 &  4.20992e-07 &     0.952243 &            1 &            1 \\
  400 &  5.23913e-09 &      0.95345 &            1 &            1 \\
  500 &  6.71114e-11 &     0.954231 &            1 &            1 \\
  600 &  8.75867e-13 &     0.954782 &            1 &            1 \\
  700 &  1.15812e-14 &     0.955193 &            1 &            1 \\
  800 &  1.54618e-16 &     0.955513 &            1 &            1 \\
  900 &  2.04046e-18 &      0.95575 &            1 &            1 \\
 1000 &  2.78188e-20 &      0.95597 &            1 &            1 \\
\end{tabular}
\end{center}

 In the next table, we treat the  ``stabilization" of the previous
automorphism. We see that the curl of the  ``stabilization"
is apparently smaller.

$$
\left\{
\begin{array}{l}
x \rightarrow xy \\
y \rightarrow y \\
z \rightarrow z
\end{array}
\right.
$$

\begin{center}
\begin{tabular}[]{r||l|l|l|l}
 \mbox{ n} & \mbox{\tiny CURL\_RATIO} & \mbox{\tiny CURL\_ROOT} & \mbox{\tiny FLUX\_RATIO} & \mbox{\tiny FLUX\_ROOT} \\
 \hline
   10 &     0.220658 &      0.85975 &     0.779342 &     0.975378 \\
   20 &    0.0832884 &     0.883139 &     0.916712 &     0.995661 \\
   50 &   0.00616004 &     0.903216 &      0.99384 &     0.999876 \\
  100 &  0.000106955 &     0.912624 &     0.999893 &     0.999999 \\
  200 &  4.26719e-08 &     0.918651 &            1 &            1 \\
  300 &  1.93205e-11 &     0.921057 &            1 &            1 \\
  400 &  9.23441e-15 &     0.922388 &            1 &            1 \\
  500 &  4.52035e-18 &     0.923231 &            1 &            1 \\

\end{tabular}
\end{center}

 In the next table, we treat the  square of the first automorphism.
 We see that the curl of the  square is apparently smaller
than that of the automorphism itself.

$$
\left\{
\begin{array}{l}
x \rightarrow xy^2 \\
y \rightarrow y
\end{array}
\right.
$$

\begin{center}
\begin{tabular}[]{r||l|l|l|l}
 \mbox{ n} & \mbox{\tiny CURL\_RATIO} & \mbox{\tiny CURL\_ROOT} & \mbox{\tiny FLUX\_RATIO} & \mbox{\tiny FLUX\_ROOT} \\
 \hline
   10 &     0.143331 &     0.823444 &     0.856670 &     0.984649 \\
   20 &    0.0408009 &     0.852184 &     0.959199 &     0.997919 \\
   50 &   0.00133009 &     0.875947 &      0.99867 &     0.999973 \\
  100 &  5.98358e-06 &     0.886686 &     0.999994 &            1 \\
  200 &  1.61895e-10 &       0.8934 &            1 &            1 \\
  300 &  4.98636e-15 &     0.896037 &            1 &            1 \\
  400 &  1.36942e-19 &     0.897101 &            1 &            1 \\
\end{tabular}
\end{center}

\bigskip

\noindent {\bf Acknowledgement}
\medskip

 We are  indebted to Alexander Ushakov
for providing us with the experimental data for  Section \ref{estim}.
His engineous algorithm allowed him to compute the curl function
$Curl_\varphi(n)$ of $\varphi: x \rightarrow xy,
~y \rightarrow y$ for values of  $n$ up to 1000, which is quite remarkable
given that in the group $F_2$,   ~$|B_{1000}|> 3^{1000} > 10^{200}$.

\bigskip

\noindent {\it http://www.math.mcgill.ca/\~\/alexeim/}
\medskip

\noindent {\it http://www.sci.ccny.cuny.edu/\~\/shpil/ \\

\end{document}